\documentstyle[11pt]{article}

\newcommand{\bbC}{{\mathbf C}}

\newcommand{\qed}{\hfill$\Box$}
\newcommand{\la}{{\lambda}} 
\newcommand{\bla}{{\bar{\lambda}}} 
\newcommand{\tr}{\mbox{\rm tr\,}} 
\newcommand{\trt}{{\epsilon(\bla)}} 
\newcommand{\too}{{\longrightarrow}}
\newcommand{\mmod}{\mbox{\rm mod\ }} 
\newcommand{\rht}{{\sl RHT}} 
\newcommand{\sign}{{\sl sign}} 
\newcommand{\core}{{\sl core}} 
\newcommand{\quot}{{\sl quot}} 
\newcommand{\inv}{{\sl inv}} 
 
\newcommand{\Aut}{{\rm Aut}}
\newcommand{\htau}{{\hat\tau}}
\newcommand{\si}{{\sigma}}
\newcommand{\bsi}{{\bar\si}} 
\newcommand{\tsi}{{\tilde\si}} 
\newcommand{\lala}{{(\la_1,\ldots,\la_k)}}
\newcommand{\mumu}{{\mu_1,\ldots,\mu_k}}
\newcommand{\TT}{{T_1,\ldots,T_k}}

\newcommand{\hgt}{{\sl ht}}
\newcommand{\Utm}{{U^{\otimes m}}}

\newlength{\zw}
\newsavebox{\obox}
\savebox{\obox}[\zw][c]{$1$}
\newsavebox{\zbox}
\savebox{\zbox}[\zw][c]{$0$}
\newsavebox{\bobox}
\savebox{\bobox}[\zw][c]{$\hat{1}$}
\newsavebox{\bzbox}
\savebox{\bzbox}[\zw][c]{$\hat{0}$}
\newsavebox{\emptybox}
\savebox{\emptybox}[\zw][c]{$\,$}
\newsavebox{\medbox}
\savebox{\medbox}[\zw][c]{$|$}
\newcommand{\uo}{\usebox{\obox}}
\newcommand{\uz}{\usebox{\zbox}}
\newcommand{\bo}{\usebox{\bobox}}
\newcommand{\bz}{\usebox{\bzbox}}
\newcommand{\ue}{\usebox{\emptybox}}

\newtheorem{thm}{Theorem}[section]
\newtheorem{pro}[thm]{Proposition}
\newtheorem{lem}[thm]{Lemma}

\newtheorem{cor}[thm]{Corollary}

\newtheorem{exa}[thm]{Example}

\newenvironment{proof}{\smallskip\par\noindent{\bf Proof. }}{\par\qed\bigskip}

\begin{document}

\title{Rim Hook Tableaux and Kostant's $\eta$-Function Coefficients}
\bibliographystyle{acm}
\author{Ron M. Adin%
\thanks{Department of Mathematics and Statistics, Bar-Ilan University,
Ramat-Gan 52900, Israel. 
Email: {\tt radin@math.biu.ac.il} } 
\thanks{Research supported in part by the Israel Science Foundation, founded by 
the Israel Academy of Sciences and Humanities, and by an internal research grant 
from Bar-Ilan University.}
\and Avital Frumkin%
\thanks{Sackler School of Mathematical Sciences, Tel-Aviv University,
Ramat-Aviv, Tel-Aviv 69978, Israel. 
Email: {\tt frumkin@math.tau.ac.il} }}
\date{March 31, 2004}

\maketitle

\begin{abstract}
Using a 0/1 encoding of Young diagrams and its consequences for rim hook tableaux, 
we prove a reduction formula of Littlewood for arbitrary characters of the symmetric 
group, evaluated at elements with all cycle lengths divisible by a given integer.  
As an application, we find explicitly the coefficients in a formula of Kostant for certain 
powers of the Dedekind $\eta$-function, avoiding most of the original machinery.
\end{abstract}

\section{Introduction}

A well known classical result due to Kostant~\cite{Kostant}, based on previous work of 
Macdonald~\cite{Macdonald}, expresses certain powers of the Dedekind $\eta$-function 
as power series, summing over the irreducible representations of suitable Lie groups.  
An important factor in these expressions
was shown to obtain only 
the values $0$, $1$, and $-1$ when the group is simply laced. 
The proof used the representation theory of Lie groups.
Stated explicitly 
(for explanation of notations see Subsection~\ref{s.kostant} below),
the result is

\begin{thm}\label{t.kostant.intro}{\rm\cite[Theorem 1]{Kostant}}
For any simple, simply connected and simply laced compact Lie group $K$,
$$
\phi(x)^{\dim K} = \sum_{\bla\in D} \trt \cdot \dim V_\bla \cdot x^{c(\bla)}
$$
and also
$$
\trt\in\{0,1,-1\}\qquad(\forall\bla\in D).
$$
\end{thm}

\medskip

In this paper we give an explicit computation of $\trt$ for $K=SU(k)$ 
(i.e., Lie type $A$), leading to a very clear description of those $\bla\in D$
for which $\trt\ne 0$. Explicitly, our main result is
\begin{thm}{\rm (see Corollary~\ref{t.bla.final})}\label{t.bla.intro}
Let $\bla = (\bla_1,\ldots,\bla_{k-1})$ be a dominant $SU(k)$-weight, and define
$\bla_k := 0$. If the numbers $(\bla_i + k - i)_{i=1}^{k}$ 
have $k$ distinct residues $(\mmod k)$
then $\trt$ is the sign of the permutation needed to transform the sequence 
of residues $(\mmod k)$ of $(\bla_i -\bar{p} + k - i)_{i=1}^{k}$ into 
the sequence $(k-i)_{i=1}^{k}$, where $\bar{p} := (\bla_1 + \ldots + \bla_k)/k$.
In all other cases, $\trt = 0$.
\end{thm}

The proof uses
combinatorial tools related to the symmetric group.  
The connection is made via the Schur-Weyl double commutant theorem.
Our arguments use the 0/1 encoding of partitions,
a corresponding interpretation of rim hook tableaux, 
and a result (known to Littlewood) concerning the computation of $S_n$-character 
values on permutations with all cycle lengths divisible by a given integer. 

The structure of the paper is as follows.
Section~\ref{s.prelim} contains some preliminaries.
In Section~\ref{s.01} we describe the 0/1 encoding of partitions and the concepts of 
$k$-quotient and $k$-core.
Section~\ref{s.decomp} deals with rim hook tableaux; the zero permutation and the
$k$-quotient of a rim hook tableau are defined and related to each other.
Finally, in Section~\ref{s.proof} we compute Kostant's coefficients for type $A$.

\section{Preliminaries}\label{s.prelim}

In this section we collect  some definitions,
notation and results that will be used in the rest of this paper.

\subsection{Partitions and Characters}\label{s.prelim.part}

A {\em partition} is a weakly decreasing infinite sequence
$\la = (\la_1,\la_2,\ldots)$ of nonnegative integers
$\la_1 \ge \la_2 \ge \ldots \ge 0$
in which only finitely many terms are nonzero.
Its {\em size} is
$$|\la| := \sum_i \la_i$$
and its {\em length} is
$$\ell(\la) := \max\,\{i\ge 1\,|\,\la_i>0\}.$$
By definition, the {\em empty partition} $\emptyset=(0,0,\ldots)$ has length 
(and size) zero.
If $|\la| = n$ then we say that $\la$ is a partition of $n$, 
and write $\la \vdash n$.

A partition $\mu$ with $m_1$ ``1''s, $m_2$ ``2''s, etc., can also be written in
multiset notation:  $\mu = (\ldots 2^{m_2} 1^{m_1})$.
For such a partition, let
\[
z_{\mu} := {m_1}!\,1^{m_1} \cdot {m_2}!\,2^{m_2} \cdots.
\]
Naturally, $z_{\mu}:=1$ for the empty partition $\mu=\emptyset$.

The {\em cycle type} of a permutation $\pi\in S_n$ is the multiset of sizes of 
cycles in $\pi$.  We can (and will) order this multiset into a partition of $n$.

\begin{pro}\label{t.cycletype}
If $\mu$ is a partition of $n$ then the number of permutations in $S_n$ having
cycle-type $\mu$ is $n\,!/z_{\mu}$.
\end{pro}

Let $\chi^{\la}$ be the irreducible character of $S_n$ indexed by $\la\vdash n$.
For $\mu\vdash n$, let $\chi_{\mu}^{\la}$ be the value of that character on an element
of $S_n$ having cycle type $\mu$.  The orthogonality property of characters, together
with Proposition~\ref{t.cycletype}, yield

\begin{pro}\label{t.charsum}
For any partition $\la$ of $n$ ($n\ge 0$),
$$
\sum_{\mu\vdash n} \frac{1}{z_{\mu}} \chi_{\mu}^{\la} =
\frac{1}{n!} \sum_{\si\in S_n} \chi^{\la}(\si) = \delta_{\la,(n)},
$$
where $\delta$ is the Kronecker delta and $(n)$ is the partition of $n$ with 
(at most) one part.
\end{pro}

\subsection{Weak Compositions and Rim Hook Tableaux}\label{s.prelim.rht}

The definitions below follow~\cite[pp.\ 345--347]{StaEC2}.

A {\em weak composition} of a nonnegative integer $n$ is an infinite sequence $\mu$, 
not necessarily decreasing, of nonnegative integers adding up to $n$. 
We write $|\mu| = n$.
If $\mu=(\mu_1,\mu_2,\ldots)$ is a weak composition and $k$ is a positive 
integer, let $k\mu := (k\mu_1,k\mu_2,\ldots)$.

A {\em rim hook} (or {\em border strip}, or {\em ribbon}) 
is a connected skew shape containing no $2 \times 2$ square.
The {\em length} of a rim hook is the number of boxes in it, and its 
{\em height} is one less than its number of rows.  By convention,
the height of an empty rim hook is zero.
 
Let $\la$ be a partition and $\mu = (\mu_1, \mu_2, \ldots)$ be a weak composition
of the same integer $n$.
A {\em rim hook tableau} of {\em shape} $\la$ and {\em type} $\mu$ is an assignment 
of positive integers to the boxes of the Young diagram of $\la$ such that
\begin{enumerate}
\item
every row and column is weakly increasing;
\item
each positive integer $i$ appears $\mu_i$ times; and
\item
the set of boxes occupied by $i$ forms a (possibly empty) rim hook.
\end{enumerate}

Equivalently, we may think of a rim hook tableau of shape $\la$ as a sequence 
$\emptyset=\la^0 \subseteq \la^1 \subseteq \cdots \subseteq \la^r = \la$
of partitions such that each skew shape $\la^i / \la^{i-1}$ is a rim hook of
length $\mu_i$ (including the empty rim hook $\emptyset$, when $\mu_i=0$). 
The {\em height} $\hgt(T)$ of a rim hook tableau $T$ is the sum of heights of all
(nonempty) rim hooks in it. 

Let $\rht_{\mu}^{\la}$ be the set of all rim hook tableaux of shape $\la$ and
type $\mu$.  For example,
$$
T= \begin{array}{ccc}
1 & 1 & 4\\
3 & 4 & 4\\
3 &   &  \\
\end{array}
\in \rht_{\mu}^{\la}
$$
where $\la=(3,3,1)$ and $\mu=(2,0,2,3)$.

\begin{pro}{\bf (Murnaghan-Nakayama Formula)}\label{t.mnf} {\rm\cite[p.~60]{JK}}
For any partition $\la$ and weak composition $\tilde\mu$ of the same integer,
$$
\sum_{T\in \rht_{\tilde\mu}^{\la}} (-1)^{\hgt(T)} = \chi_{\mu}^{\la}
$$
where $\mu$ is the unique partition obtainable from $\tilde\mu$ by reordering 
the parts.
\end{pro}

It follows that the sum on the left-hand side is invariant under reordering the
weak composition $\tilde\mu$.  The size of the set $\rht_{\tilde\mu}^{\la}$ 
itself is, however, not invariant; this set may even be empty for some 
orderings and nonempty for others.

\section{The 0/1 Encoding of Partitions}\label{s.01}

In this section we shall describe an encoding of partitions (or Young diagrams) 
in terms of 0/1 sequences, as well as the concepts of $k$-quotient and $k$-core
of a partition.  
This will lead, in Section~\ref{s.decomp}, to a very transparent description of 
rim hook tableaux.
This description, in turn, will be used in Section~\ref{s.proof} to give a reduction of 
a character sum over a coset in a large group to a similar sum in a smaller group, 
via the Murnaghan-Nakayama character formula.
For historical and other remarks regarding the 0/1 encoding,
see Subsection~\ref{s.01remarks} below.

\subsection{The 0/1 Sequence of a Partition}\label{s.01.sub}

Let $\la=(\la_1,\la_2,\ldots)$ be a partition of length $n$, so that 
$\la_1\ge\ldots\ge\la_n > 0 = \la_{n+1} = \ldots$ $(n\ge 0)$.  
This is usually denoted graphically by a {\em Young diagram}, which has $n$
left-justified rows of boxes, with $\la_i$ boxes in row $i$ ($1\le i\le n$).  
Extend the South-Eastern borderline of the diagram by a vertical ray extending
down (South) from the South-Western corner, and a horizontal ray extending to the 
right (East) from the North-Eastern corner.  Denote the South-Eastern borderline, 
together with these two rays, by $\partial\la$, and call it the {\em extended borderline} 
of $\la$.

Now trace the edges of $\partial\la$, starting from the vertical ray and
ending at the horizontal ray.  Each step is either up (North) or to the right (East).
Encode this walk by a 0/1 sequence $s(\la)$, where each ``up" step is encoded by ``0"
and each ``right" step is encoded by ``1".  The resulting doubly-infinite sequence
is initially all ``0"s (corresponding to the vertical ray) and eventually all ``1"s
(corresponding to the horizontal ray).

Each digit (``0" or ``1") in the sequence corresponds to an edge of the extended
borderline $\partial\la$, and each ``space" between digits corresponds to a vertex 
of $\partial\la$.  A special ``space" in the sequence corresponds to the unique vertex
of $\partial\la$ which is on the same diagonal as the North-Western corner of $\la$.  
Call this space the {\em median} of $s(\la)$ and denote it by ``$|$"; it is the unique space
satisfying
\begin{equation}\label{e.median}
\#\{\hbox{\rm ``1''s preceding the median}\} = 
\#\{\hbox{\rm ``0''s succeeding the median}\}.
\end{equation}

Many natural parameters of a partition $\la$ can be easily described using $s(\la)$.
The ``0''s of $s(\la)$ correspond (in reverse order) to the rows of the diagram of $\la$. 
The length $\la_i$ of the $i$th row is the number of ``1''s to the left of the $i$th (from the
right) ``0'' in $s(\la)$.  In particular, the empty rows correspond to the infinite sequence
of leading ``0''s in $s(\la)$.  
Similarly, the length $\la'_j$ of the $j$th column is the number of ``0''s 
to the right of the $j$th (from the left) ``1'' in $s(\la)$.
The length $\ell(\la)=\la'_1$ of $\la$ is the number of ``0''s to the right of the first ``1'' 
in $s(\la)$.
The size of $\la$ is the number of {\em inversions} in $s(\la)$:
\begin{equation}\label{e.size}
|\la| = \#\{(i,j) \,|\, i<j,\, s_i=1,\, s_j=0\}.
\end{equation}
The side-length $r$ of the Durfee square of $\la$ (the largest square contained in the 
diagram of $\la$) is the value of either side of equation~(\ref{e.median}).
The Frobenius notation $\la=(\alpha_1,\ldots,\alpha_r \,|\, \beta_1,\ldots,\beta_r)$
is obtained by letting $\alpha_i$ (respectively, $\beta_i$) be the distance from the median
to the $i$th from the right (left) ``0'' (``1''); only ``0''s to the right (left) of the median are 
taken into account here.
The empty partition corresponds, of course, to the sequence $\ldots\uz\uz | \uo\uo\ldots$.

\subsection{The $k$-Quotient}\label{s.quotient}

For an integer $k\ge 2$, define the {\em $k$-quotient} of a partition $\la$ by
\[
\quot_k(\la) := (\la_1,\ldots,\la_k),
\]
where $\la_i$ $(1\le i\le k)$ is the partition whose 0/1 sequence $s(\la_i)$ is the 
subsequence containing every $k$th digit of $s(\la)$, including the $i$th digit after 
the median of $s(\la)$.

\begin{exa} 
Let $\la=(4,4,3,2)$ and $k=3$.  Then $\quot_3(\la) = (\la_1,\la_2,\la_3)$ where
$\la_1=(1,1)$, $\la_2=\emptyset$ and $\la_3=(2)$, since  
$$
s(\la) = \ldots\uz\uz\uo\uo\uz\uo | \uz\uo\uz\uz\uo\uo \ldots \;\too\;
\begin{array}{c}
\ldots\uz\ue\ue\uo\ue\ue | \uz\ue\ue\uz\ue\ue \ldots = s(\la_1)\\
\ldots\ue\uz\ue\ue\uz\ue | \ue\uo\ue\ue\uo\ue \ldots = s(\la_2)\\
\ldots\ue\ue\uo\ue\ue\uo | \ue\ue\uz\ue\ue\uo \ldots = s(\la_3)\\
\end{array}
$$
\end{exa}

Note that the medians shown in the example are the ones ``inherited'' from 
$s(\la)$, and do not necessarily coincide with the intrinsically defined medians 
of the various $s(\la_i)$.  
Actually, in order to reconstruct $\la$ from its $k$-quotient we need to know 
the relative positions of the medians of the various subsequences $s(\la_i)$.  
This information is provided by the $k$-core, to be defined in the sequel.

\subsection{Removing a Rim Hook}

Removing a rim hook of length $k$ from a Young diagram $\la$ is equivalent to finding
a ``1" and a ``0" at distance $k$ apart in $s(\la)$, with the ``1" preceding the ``0",
and interchanging these two digits.  The height $h$ of the rim hook is equal to 
the number of ``0''s strictly between the two digits to be interchanged.
\begin{exa}
Denoting by ``$\,\bo$'' and ``$\,\bz$'' the two digits to be interchanged,
 here is an example of removing a rim hook of length $4$ and height $1$:
$$
\ldots \uz\uz\bo\uo\uz\uo | \bz\uo\uz\uz\uo\uo \ldots \,\too\, 
\ldots \uz\uz\bz\uo\uz\uo | \bo\uo\uz\uz\uo\uo \ldots  
$$
\end{exa}

It follows that removing from $\la$ a rim hook of length divisible by $k$ affects
only one of the subsequences $s(\la_1),\ldots,s(\la_k)$, where
$(\la_1,\ldots,\la_k) = \quot_k(\la)$.
In particular, removing a rim hook of length exactly $k$ is equivalent to 
a transposition $10 \to 01$ of {\em consecutive} digits in one of these 
subsequences.  Note also that, by condition~(\ref{e.median}), this operation 
changes neither the position of the median of $s(\la)$ nor the intrinsic median of 
any subsequence $s(\la_i)$.

\subsection{The $k$-Core}

The {\em $k$-core} of a partition $\la$, $\core_k(\la)$, is the unique partition 
obtained from $\la$ by removing as many rim hooks of length $k$ as possible.  
By the remarks in the previous subsection, this corresponds to ``squeezing'' 
all the ``0''s to the left  (and all the ``1''s to the right) 
in each of the subsequences $s(\la_1),\ldots,s(\la_k)$ separately.  
This interpretation shows that the $k$-core is well-defined, i.e., does not depend 
on the choice of rim hooks to be removed.  It also follows that the $k$-core has 
trivial $k$-quotient, i.e.,
$$
\quot_k(\core_k(\la)) = (\emptyset,\ldots,\emptyset)\qquad(\forall\la).
$$
Thus, the $k$-core encodes only the relative positions of the intrinsic medians 
of the various subsequences $(\mmod k)$ of $s(\la)$; their ``average'' is the
median of $s(\la)$.  By the comment at the end of Subsection~\ref{s.quotient}
we now get the following result.

\begin{lem}{\rm\cite[Theorem 2.7.30]{JK}}
For any $k\ge 2$, a partition is uniquely determined by its $k$-core and $k$-quotient. 
\end{lem}

The following consequence of the remarks above is noted here for future reference.

\begin{lem}{\bf (Empty $k$-Core Criterion)} \label{t.empty.core}\\
A partition $\la$ has an empty $k$-core if and only if the median of 
each subsequence $(\mmod k)$ of $s(\la)$ coincides with the median of $s(\la)$.
\end{lem}

\begin{exa}\rm $(k=2)$\\
Each of the following two examples shows the maximal rim hook removal 
from a 0/1 sequence and the corresponding operation on its two subsequences 
$(\mmod 2)$.
\begin{enumerate}
\item
Empty $2$-core:
$$
\begin{array}{ccc}
 \ldots \uo\uz\uo\uo | \uz\uz\uz\uo \ldots &\too
&\ldots \uz\uz\uz\uz | \uo\uo\uo\uo \ldots\\
&&\\
 \ldots \uo\ue\uo\ue | \uz\ue\uz\ue \ldots &\too
&\ldots \uz\ue\uz\ue | \uo\ue\uo\ue \ldots\\
 \ldots \ue\uz\ue\uo | \ue\uz\ue\uo \ldots &\too
&\ldots \ue\uz\ue\uz | \ue\uo\ue\uo \ldots\\
\end{array}
$$
\item
Nonempty $2$-core:
$$
\begin{array}{ccc}
 \ldots \uo\uz\uo\uo | \uz\uz\uo\uz \ldots &\too
&\ldots \uz\uz\uo\uz | \uo\uz\uo\uo \ldots\\
&&\\
 \ldots \uo\ue\uo\ue | \uz\ue\uo\ue \ldots &\too
&\ldots \uz\ue\uo\ue | \uo\ue\uo\ue \ldots\\
 \ldots \ue\uz\ue\uo | \ue\uz\ue\uz \ldots &\too
&\ldots \ue\uz\ue\uz | \ue\uz\ue\uo \ldots\\
\end{array}
$$
\end{enumerate}
\end{exa}

We shall also need the following quantitative relation.
\begin{lem}\label{t.size}
If 
$$
\quot_k(\la) = (\la_1,\ldots,\la_k)
$$
and
$$
\core_k(\la) = \la_0
$$
then 
$$
|\la| = |\la_0| + k \cdot (|\la_1| + \ldots + |\la_k|).
$$
\end{lem}
\begin{proof}
Recall expression~(\ref{e.size}) for the size of a partition in terms of 
the number of inversions in its 0/1 sequence.  Removal of a rim hook 
of length $k$ from $\la$ corresponds to an adjacent  transposition $10 \to 01$
in one of the subsequences $s(\la_i)$, thus reducing by $1$ its inversion
number (and the size of $\la_i$), not affecting any of the other subsequences.  
This removal clearly reduces by $k$ the size of $\la$, as can also be seen using $s(\la)$.
\end{proof}

Note the following special case, which actually follows directly from the 
definition of $k$-core.

\begin{cor}
If $\la$ is a partition with $\core_k(\la) = \emptyset$ then $k$ divides $|\la|$.
\end{cor}

\subsection{Historical Remarks}\label{s.01remarks}

The 0/1 encoding is a succinct and useful way to represent a Young diagram.
It was rediscovered many times; 
e.g., by Shahar Mozes and the present authors in the 1980's (unpublished). 
See \cite[pp.~467, 517]{StaEC2} for references, dating from 1959 on. 
An essentially equivalent technique, using bead configurations and abaci, appears 
in~\cite[Ch.~2.7]{JK}.  
Fairy sequences (or, equivalently, content vectors), used for example in~\cite{FS}, 
are essentially partial sums of 0/1 sequences.
The introduction of the median and the statement of Lemma~\ref{t.empty.core} 
(the empty $k$-core criterion) seem to be new.

\section{The Decomposition Theorem}\label{s.decomp}

\subsection{The $k$-Quotient of a Rim Hook Tableau}\label{s.quot.tab}

Using the 0/1 encoding of partitions, a rim hook tableau $T \in \rht_\mu^\la$ 
corresponds to a finite sequence of operations 
of either the form $1 \ldots 0 \to 0 \ldots 1$ or the form ``do nothing'', 
transforming $s(\la)$ into $s(\emptyset) = \ldots 00|11 \ldots$, and
whose sequence of lengths is the weak composition $\mu$ in reverse order.
Note that we view $T$ here as ``peeling'' rather than ``constructing'' the diagram 
$\la$, and that empty rim hooks are allowed.
If $\mu$ is a partition (i.e., weakly decreasing) then rim hooks are removed in order 
of weakly increasing length.

Let $\la$ be a partition with empty $k$-core, and let 
$\quot_k(\la) = (\la_1,\ldots,\la_k)$ be its $k$-quotient.
If $|\la_i| = p_i$ $(1\le i\le k)$ then, by Lemma~\ref{t.size}, $p_1 + \ldots + p_k = p$
where $|\la| = kp$.  
Let $T$ be a rim hook tableau of shape $\la$, with all rim hook lengths divisible by $k$.  
Its type has the form $k\mu$, with $\mu$ a weak composition of $p$.
In the sequel we shall usually assume, for simplicity, that $\mu$ is actually 
a partition.
Clearly, for any $t\ge 0$, removing a rim hook of length $kt$ from $\la$ is equivalent to 
removing a rim hook of length $t$ from a suitable $\la_i$.  
The peeling of $T$ thus corresponds to peeling $k$ rim hook tableaux $\TT$.
If $\mu$ is a partition then the type of each $T_i$ is also a partition.
We can then define a mapping 
$$
\omega_{k\mu}^{\la} : \rht_{k\mu}^{\la} \too \bigcup_{(\mumu)\in P_{\mu}^{\la}} 
\left(\rht_{\mu_1}^{\la_1} \times \cdots \times \rht_{\mu_k}^{\la_k} \right)
$$
where $P_{\mu}^{\la}$ is the set of all $k$-tuples $(\mumu)$ of partitions satisfying
$$
|\mu_i| = |\la_i| \qquad(1\le i\le k)
$$
and
$$
\mu_1  \oplus \ldots \oplus \mu_k = \mu.
$$
Here the direct sum of $\mumu$ is the partition obtained by reordering (in 
weakly decreasing order) all the parts of all the $\mu_i$.

Recall the definition of $z_{\mu}$ from Subsection~\ref{s.prelim.part}.

\begin{lem}\label{t.rht.dec}
If $\la$ and $\mu$ are partitions satisfying $\core_k(\la)=\emptyset$ and
$|\la| = k|\mu|$ then the map $\omega_{k\mu}^{\la}$ defined above is surjective, 
and each $k$-tuple  
$$
(\TT) \in \rht_{\mu_1}^{\la_1} \times \cdots \times \rht_{\mu_k}^{\la_k}
$$
in its range is obtained from exactly 
$z_{\mu}/(z_{\mu_1} \cdots z_{\mu_k})$ elements (rim hook tableaux) in its domain. 
\end{lem}
\begin{proof}
Let $\mu = (\ldots 2^{m_2} 1^{m_1})$.  Given $(\TT)$ we can reconstruct $T$ by 
successively adding the rim hooks of all the $T_i$.  The only ambiguity is in
deciding, for each $j\ge 1$, which of the $m_j$ rim hooks of length $j$ in $T$
should come from each $T_i$ ($1\le i\le k$).  The number of possibilities is
a product (over $j$) of multinomial coefficients, and boils down to the claimed
formula.
\end{proof}

\subsection{The Zero Permutation}

Let $\la$ be a partition.  Label the zeros in $s(\la)$ by distinct labels, say by the
positive integers. It will actually suffice to label only the finitely many 
zeros corresponding to non-empty rows of $\la$, i.e., the zeros succeeding the 
first ``1'' in $s(\la)$.  Let $\ell(\la)$ be the length of $\la$, and 
let $T$ be a rim hook tableau of shape $\la$ and arbitrary type $\mu$. 
The {\em zero permutation} of $T$ is the permutation $\pi_T\in S_{\ell(\la)}$ 
of the labeled zeros in the sequence $s(\emptyset)$ obtained from $s(\la)$ by 
successively removing rim hooks according to $T$, namely in order of 
decreasing entries.

\begin{exa}
Let
$$
T= \begin{array}{ccc}
1 & 1 & 4 \\
3 & 4 & 4 \\
3 &   &   \\
\end{array}
\in \rht_{\mu}^{\la}
$$
where $\la=(3,3,1)$ and $\mu=(2,0,2,3)$.
Label the zeros in $s(\la)$, and remove rim hooks from $\la$ according to $T$
(starting with the rim hook with entries equal to 4):
\begin{eqnarray*}
\ldots 1 0_1 \hat{1} | 1 0_2 \hat{0}_3 \ldots &\stackrel{4}{\too}&
\ldots \hat{1} 0_1 \hat{0}_3 | 1 0_2 1 \ldots\stackrel{3}{\too}
\ldots 0_3 0_1 1 | 1 0_2 1\\
&\stackrel{2}{\too}&
\ldots 0_3 0_1 \hat{1} | 1 \hat{0}_2 1 \ldots\stackrel{1}{\too}
\ldots 0_3 0_1 0_2 | 1 1 1 \ldots;
\end{eqnarray*}
hence $\pi_T = 312 \in S_3$.
\end{exa}

Recall the definition of $\hgt(T)$ from Subsection~\ref{s.prelim.rht}.

\begin{lem}\label{t.sign.hgt}
The sign of the zero permutation
$$
\sign(\pi_T) = (-1)^{\hgt(T)}.
$$
\end{lem}

\begin{proof}
Recall that 
$$
\sign(\pi_T) = (-1)^{\inv(\pi_T)},
$$
where the {\em inversion number} of a permutation $\pi$ is
$$
\inv(\pi) := \#\{(i,j)\,|\,i<j,\,\pi(i)>\pi(j)\}.
$$

Both $\inv(\pi_T)$ and $\hgt(T)$ are equal to zero for the empty tableau 
$T=\emptyset$.  They both change by $h (\mmod 2)$ when we remove from 
a tableau $T$ a nonempty rim hook of length $k$ and height $0\le h\le k-1$
(since $\pi_T$ records only the zeros of $s(\la)$, and during the rim hook removal 
one of these zeros ``skips'' exactly $h$ others). 
Therefore
$$
\inv(\pi_T) \equiv \hgt(T) \quad(\mmod 2)
$$
and the proof is complete.
\end{proof}

\begin{cor}\label{t.sign_SYT}
If $T$ is a rim hook tableau of type $(1^n)$ (i.e., a standard Young tableau)
then $\sign(\pi_T) = 1$.
\end{cor}

\begin{proof}
The height of a rim hook of length $1$ is $0$, so that $\hgt(T)=0$ here.
Indeed, the permutation $\pi_T$ is the identity permutation, since 
removing a rim hook of length $1$ is an adjacent transposition $10 \to 01$ 
in $s(\la)$ (where $\la$ is the shape of $T$), so does not affect the order
of zeros.
\end{proof}

\subsection{The Sign Decomposition Theorem}

This subsection is devoted to the proof of Theorem~\ref{t.sign} below.
An essentially equivalent formulation of this theorem, in terms of 
character values, 
was given by Littlewood~\cite[pp.~143--146]{L}.

We first state an important special case as a lemma.

\begin{lem}\label{t.zero_perm}
Let $\la$ be a partition with empty $k$-core, and let $p := |\la|/k$.
Then, for any $T\in \rht_{(k^p)}^{\la}$, the zero permutation $\pi_T$ depends
only on $\la$ (and not on the choice of $T$), and will be denoted $\pi_{\la}$.
In particular, 
$$
\chi_{(k^p)}^{\la} = \sign(\pi_{\la}) \cdot |\rht_{(k^p)}^{\la}|.
$$
\end{lem}
\begin{proof}
$T$ is of shape $\la$ and type $(k^p)$.
Removing a rim hook of length $k$ from $s(\la)$ corresponds to 
an adjacent transposition $10 \to 01$ in a suitable subsequence $s(\la_i)$, 
and does not change the zero permutation in this subsequence.  
Thus, the effect of removing all the rim hooks in $T$ is merely to shift, in each 
subsequence $s(\la_i)$, all the zeros as much to the left as possible, 
without changing their order.  The resulting zero permutation $\pi_T$, therefore, 
depends only on the partition $\la$ and not on the specific choice of $T$.
Denote this permutation by $\pi_{\la}$.
The claimed formula follows from Proposition~\ref{t.mnf} 
(with $\mu = \tilde\mu = (k^p)$), Lemma~\ref{t.sign.hgt} and the definition 
of $\pi_{\la}$.
\end{proof}

\begin{thm}\label{t.sign} {\bf (Sign Decomposition Theorem)\\}
Let $\la$ be a partition with empty $k$-core, $\mu$ a partition of 
$p := |\la|/k$, and $T\in \rht_{k\mu}^{\la}$.
If $\omega_{k\mu}^{\la}(T) = (\TT)$ then
$$
\sign(\pi_T) = \sign(\pi_{\la}) \cdot \sign(\pi_{T_1}) \cdots \sign(\pi_{T_k}),
$$
where $\pi_{\la}$ is as defined in Lemma~\ref{t.zero_perm}.
\end{thm}

In our proof of this theorem we shall use the following result of 
Garsia and Stanley~\cite{Sta84LAMA}. We note that this result may also be 
proved using 0/1 sequences.

\begin{lem}{\rm\cite[Lemma 7.3]{Sta84LAMA}}\label{t.gasta} 
If $H$ is a skew shape which is a rim hook of length $kt$, then there exists a skew 
rim hook tableau of the shape of $H$ consisting of $t$ rim hooks of length $k$ each.
\end{lem}

\begin{cor}
If $\rht_{k\mu}^{\la} \ne \emptyset$ for some weak composition $\mu$ of size 
$p := |\la|/k$, then this also holds for the partition $\mu = (1^p)$.
\end{cor}

\smallskip\par\noindent{\bf Proof of Theorem~\ref{t.sign}. }
$T$ is a rim hook tableau of shape $\la$, with all rim hook lengths
divisible by $k$.  Consider a rim hook in $T$, of length $kt$ ($t\ge 1$).  
By Lemma~\ref{t.gasta}, we can replace this rim hook by a suitable sequence of 
$t$ rim hooks of length $k$ each. 
The resulting set of zero positions is the same, whether we remove the original 
long rim hook or the $t$ short ones; the two zero permutations, though, differ by a cycle 
of length $h+1$, where $h$ is the height of the rim hook of length $t$ corresponding, 
in one of the tableaux $T_i$, to the original rim hook of length $kt$.

\begin{exa}
$k=2,\,t=3,\,h=1.$\\
Removing one rim hook of length $kt = 6$ from $s(\la)$:
$$
\ldots \hat{1} 1 0_1 0_2 | 1 0_3 \hat{0}_4 1 \ldots\too
\ldots 0_4 1 0_1 0_2 | 1 0_3 1 1 \ldots
$$
Removing the corresponding rim hook of length $t=3$ and height $h=1$ from $s(\la_1)$:
$$
\ldots \hat{1} 0_1 | 1 \hat{0}_4 \ldots\too
\ldots 0_4 0_1 | 1 1 \ldots
$$
Removing $t=3$ rim hooks of length $k=2$ each from $s(\la)$:
\begin{eqnarray*}
\ldots \hat{1} 1 \hat{0}_1 0_2 | 1 0_3 0_4 1\ldots &\too&
\ldots 0_1 1 1 0_2 | \hat{1} 0_3 \hat{0}_4 1\ldots  \too
\ldots 0_1 1 \hat{1} 0_2 | \hat{0}_4 0_3 1 1\ldots \\
 &\too& 
\ldots 0_1 1 0_4 0_2 | 1 0_3 1 1\ldots 
\end{eqnarray*}
The two resulting zero permuations, $4123$ and $1423$, differ by a cycle of 
length $h+1=2$ (i.e., a transposition).
\end{exa}

It follows that this replacement multiplies both $\sign(\pi_T)$ and the corresponding
$\sign(\pi_{T_i})$ by $(-1)^h$, without changing $\sign(\pi_{\la})$.  Iterating
this process, we eventually get to the case of $T$ with all rim hooks of the same
length $k$, and Lemma~\ref{t.zero_perm} completes the proof of our claim.

\par\qed

\section{Computation of Kostant's Coefficients}\label{s.proof}

\subsection{Kostant's Formula}\label{s.kostant}

Let us recall Kostant's formula for powers of the $\phi$-function
$$
\phi(x) := \prod_{n=1}^{\infty} (1-x^n),
$$
which is a slight modification of the Dedekind $\eta$-function
$$
\eta(x) := x^{1/24} \phi(x).
$$
There are nice power-series expansions of
$\phi(x)^{-1}$ (the generating function for partition numbers),
$\phi(x)$ (due to Euler), $\phi(x)^3$ (due to Jacobi), $\phi(x)^{24}$
(whose coefficients are the values of the Ramanujan $\tau$-function), and others.
Note that $3$ and $24$ are the dimensions of the simple Lie algebras of types
$A_1$ and $A_4$ (or the simple compact Lie groups $SU(2)$ and $SU(5)$).

Kostant~\cite{Kostant}, building on previous work of Macdonald~\cite{Macdonald},
gave an algebraic interpretation to the power-series expansion of certain powers
of $\phi(x)$, including the ones mentioned above. We now state his result.

\smallskip

Let $K$ be a compact, simple, simply connected, simply laced Lie group, and let
$D$ be its set of dominant weights.
For $\bla\in D$, let $V_{\bla}$ be an irreducible $K$-module 
corresponding to $\bla$; and let
$$
c(\bla) := (\bla+\rho,\bla+\rho) - (\rho,\rho),
$$
where $\rho$ is half the sum of all the positive roots.
Let $V_\bla^T$ be the zero-weight subspace of $V_\bla$, i.e., the subspace of 
all vectors which are pointwise invariant under a fixed maximal torus $T$ in 
$K$. Let $W$ be the Weyl group of $K$, and let 
$$
\theta_\bla : W \too \Aut V_\bla^T
$$
be the representation of $W$ on $V_\bla^T$. Denote 
$$
\epsilon(\bla) := \tr \theta_\bla(\tau),
$$
where $\tau$ is any Coxeter element in $W$.
Using these notations, here is the remarkable result of Kostant:

\begin{thm}\label{t.kostant}{\rm\cite[Theorem 1]{Kostant}}
For any (simply laced) $K$ as above,
$$
\phi(x)^{\dim K} = \sum_{\bla\in D} \trt \cdot \dim V_\bla \cdot x^{c(\bla)}
$$
and also
$$
\trt\in\{0,1,-1\}\qquad(\forall\bla\in D).
$$
\end{thm}

We shall discuss in this paper the case of type $A$, i.e., $K=SU(k)$ 
(for which the Weyl group $W=S_k$, the symmetric group), and
compute explicitly the value of $\trt$ by ``elementary'' means.

\subsection{From $SU(k)$ to $S_m$}

In this section we express $\trt$ as a sum over the symmetric group $S_m$,
for suitable values of $m$.

We can replace the compact Lie group $SU(k)$ by the algebraic group 
$SL_k$, which has the same representations.  
In fact, it will be even more convenient to use the (nonsimple) group $GL_k$.
The irreducible representations of $GL_k$ are indexed by their maximal weights,
which may be represented by partitions $\la$ (of an arbitrary nonnegative 
integer) with $\ell(\la)\le k$.
The corresponding irreducible representation of $SL_k$ (or $SU(k)$) is actually
determined by the dominant weight
$\bla = (\la_1 - \la_k, \la_2 - \la_k, \ldots, \la_{k-1} - \la_k)$.
The correspondence between irreducible representations of $GL_k$ and of $SL_k$
is not one-to-one.

Let $U$ be a $k$-dimensional vector space over $\bbC$, and consider its $m$th 
tensor power $\Utm$. This tensor power carries a $GL_k$-action, via
the natural action of $GL_k$ on $U$ (after choosing a basis for $U$), and also 
an $S_m$-action $\rho$ by permuting the factors:
$$
\rho(\si)(v_1 \otimes \ldots \otimes v_m) := 
v_{\si^{-1}(1)} \otimes \ldots \otimes v_{\si^{-1}(m)}
\qquad(\forall \si\in S_m,\,v_1,\ldots,v_m \in U).
$$
The following classical result expresses Schur-Weyl duality between $S_m$ 
and $GL_k$.
\begin{pro}{\bf(Double Commutant Theorem)} {\rm\cite[p.\ 374]{GW}}\label{t.double}
$$
\Utm \cong_{S_m \times GL_k} \bigoplus_{\la\in Par_k(m)} (S^\la \otimes V_\la)
$$
is an isomorphism of $(S_m \times GL_k)$-modules,
where $\la$ runs through all partitions of $m$ with at most $k$ parts,
$S^\la$ is the irreducible $S_m$-module corresponding to $\la$, and
$V_\la$ is the corresponding irreducible $GL_k$-module.
\end{pro}

Now let $T$ be a maximal torus in $SL_k$, and $W = N(T)/T$ the corresponding 
Weyl group; thus $W \cong S_k$.
Let $(\Utm)^T$ be the zero-weight subspace of $\Utm$, 
consisting of its $T$-invariant vectors.  By definition, it carries an action
of $W=N(T)/T$.  

The following statement is clear.

\begin{lem}\label{t.basis}\
\begin{description}
\item{(a)}
$$
(\Utm)^T \ne \{0\} \iff k\,|\,m.
$$
\item{(b)}
If $m=kp$ and $\{e_1,\ldots,e_k\}$ is a basis of $U$ consisting of (simultaneous) 
eigenvectors for $T$, then
$$ 
B := \{e_{i_1} \otimes \ldots \otimes e_{i_m}\,|\,\mbox{\rm as a multiset, }
\{i_1,\ldots,i_m\} = (1^p 2^p\ldots k^p) \}
$$
is a basis of $(\Utm)^T$.  
\end{description}
\end{lem}

Now take 
$$
\htau = 
\left(\begin{array}{cccc}
0 & \cdots  & 0 & (-1)^{k-1}\\
1 & 0 &   & 0 \\
  & \ddots & \ddots & \vdots\\
  &   & 1 & 0 
\end{array}\right)
\in SL_k.
$$
Actually, $\htau\in N(T)\subseteq SL_k$ where $T$ is the maximal torus 
consisting of diagonal matrices in $SL_k$. The corresponding vectors 
$e_1,\ldots,e_k$ form the standard basis of $U \cong \bbC^k$.
The representative of $\htau$ in the Weyl group $W=N(T)/T$ corresponds to
the Coxeter element $\tau = (1,2,\ldots, k)\in S_k$ under the natural
isomorphism between the group of $k \times k$ permutation matrices and the 
symmetric group $S_k$.
Recall the definitions of
$$
\theta_\bla : W \too \Aut V_\bla^T
$$
and of $\trt$ from Subsection~\ref{s.kostant}.
Our goal is to prove the claim $\trt\in\{0,1,-1\}$ of Theorem~\ref{t.kostant} 
by ``elementary'' means, i.e., by using the 
representation theory of the symmetric group instead of highest weight theory.

For $\la\in Par_k(m)$, let $\chi^\la$ be the corresponding irreducible 
$S_m$-character. The element
\begin{equation}\label{e.ela}
e_{\la} = \frac{\chi^{\la}(id)}{m!} \sum_{\si\in S_m} \chi^{\la}(\si)\si\in \bbC S_m
\end{equation}
is a central idempotent in the group algebra \cite[p.\ 50]{Serre}.
It defines, via the representation $\rho$,  a linear operator 
$\rho(e_{\la}) \in\mbox{\rm End}_{\bbC}(\Utm)$ which is 
a projection from $\Utm$ onto its $\chi^\la$-isotypic component
\cite[p.\ 21]{Serre}, isomorphic (by Proposition~\ref{t.double})
to a direct sum of $\dim V_\la$ copies of $S^\la$. 
Note that complex conjugation of character values may (and will) be 
suppressed in this paper, since all the characters of $S_m$ are real 
(actually, integer)  valued.
Also, $\chi^{\la}(id) = \dim S^{\la}$ is the multiplicity of the 
representation $V_\la$ in $\Utm$, viewed as a $GL_k$-module.

Let $\bla$ be the dominant weight corresponding to the partition $\la$.
Since $\dim S^\la = \chi^\la(id)$,
$$
\trt
= \tr\theta_\bla(\tau) 
= \tr\left(\left.\htau\right|_{V_\la^T}\right)
= \frac{1}{\chi^\la(id)} \tr\left(\left.\htau\right|_{S^\la \otimes V_\la^T}\right).
$$
The projection $\rho(e_\la)$ acts like the identity map on 
$S^\la \otimes V_\la^T$,
and like the zero map on a complementary subspace of $(\Utm)^T$.
It follows that
\begin{equation}\label{e.trt}
\trt
= \frac{1}{\chi^\la(id)} \tr\left(\left.\htau\rho(e_\la)\right|_{(\Utm)^T}\right).
\end{equation}

Recalling Lemma~\ref{t.basis} we may assume that $k\,|\,m$, say $m=kp$.
Choose one of the elements of the basis $B$ for $(\Utm)^T$, say
\begin{equation}\label{e.v0}
v_0 := (e_1 \otimes \ldots \otimes e_1) \otimes
       (e_2 \otimes \ldots \otimes e_2) \otimes\ldots\otimes
       (e_k \otimes \ldots \otimes e_k),
\end{equation}
where each of the vectors $e_1,\ldots,e_k$ appears in $p$ consecutive factors.
Every other basis element in $B$ has the form $\rho(\si)(v_0)$ for some $\si\in S_m$.
Let the subgroup $H$ of $S_m$ be the stabilizer of $v_0$.  Clearly 
$H \cong S_p \times \ldots \times S_p$ $(k \mbox{\rm\ factors})$.  

Let $\langle\cdot,\cdot\rangle$ be the inner product on $(\Utm)^T$ for which the 
basis $B$ is orthonormal.

\begin{lem}\label{t.eps.sm}
For $\bla$ corresponding to $\la\in Par_k(m)$ as above,
$$
\trt = \frac{1}{|H|} \sum_{\si\in S_m} \chi^\la(\si) \langle \htau\rho(\si)(v_0),v_0 \rangle.
$$
\end{lem}
\begin{proof}
Let $B$ be the above basis of $(\Utm)^T$.  By~(\ref{e.trt}) and~(\ref{e.ela})
\begin{eqnarray*}
\trt 
     &=& \frac{1}{\chi^\la(id)} 
         \tr\!\left(\left.\htau\rho(e_\la)\right|_{(\Utm)^T}\right) =\\
     &=& \frac{1}{m!} \sum_{\si\in S_m} \chi^{\la}(\si)
         \,\tr\left(\left.\htau\rho(\si)\right|_{(\Utm)^T}\right) =\\
     &=& \frac{1}{m!} \sum_{\si\in S_m} \chi^{\la}(\si)
         \sum_{v\in B} \langle \htau\rho(\si)(v),v \rangle.
\end{eqnarray*}
$S_m$ acts transitively on the orthonormal basis $B$, 
with $H$ as the stabilizer of $v_0$; and $\rho(\si)$ is a 
unitary operator, for each $\si\in S_m$.  Thus
\begin{eqnarray*}
\trt
     &=& \frac{1}{m!} \sum_{\si\in S_m} \chi^{\la}(\si)
         \frac{1}{|H|} \sum_{\tsi\in S_m} 
         \langle\htau\rho(\si)\rho(\tsi)(v_0),\rho(\tsi)(v_0)\rangle =\\
     &=& \frac{1}{m!|H|} \sum_{\si\in S_m} 
         \sum_{\tsi\in S_m} \chi^{\la}(\si)
         \langle \rho({\tsi}^{-1})\htau\rho(\si\tsi)(v_0),v_0 \rangle.
\end{eqnarray*}
The actions of $\htau\in GL_k$ and $\tsi\in S_m$ commute; therefore
$$
\trt
      = \frac{1}{m!|H|} \sum_{\si\in S_m} 
         \sum_{\tsi\in S_m} \chi^{\la}(\si) 
         \langle \htau\rho({\tsi}^{-1}\si\tsi)(v_0),v_0 \rangle.
$$
Denoting 
$\si' := \tsi^{-1}\si\tsi$, we now get (summing over $\si'$ and $\tsi$ instead of 
$\si$ and $\tsi$):
\begin{eqnarray*}
\trt
     &=& \frac{1}{m!|H|} \sum_{\si'\in S_m} 
         \sum_{\tsi\in S_m} \chi^{\la}(\tsi\si'\tsi^{-1})
         \langle \htau\rho(\si')(v_0),v_0 \rangle =\\
     &=& \frac{1}{|H|} \sum_{\si'\in S_m} \chi^{\la}(\si')
         \langle \htau\rho(\si')(v_0),v_0 \rangle,
\end{eqnarray*}
as claimed.
\end{proof}

\subsection{From $S_{kp}$ to $S_p^{\times k}$}\label{s.sm.spk}

We shall now show that the summation in Lemma~\ref{t.eps.sm} is actually not on
all of $S_m = S_{kp}$, but rather on a certain coset of the subgroup defined above,
$$
H \cong S_p^{\times k} := S_p \times \ldots \times S_p\quad(k \mbox{\rm\ factors}).  
$$

\begin{lem}\label{t.sm.h}
There exists a permutation $\si_0 \in S_{kp}$ such that
$$
\trt = \frac{(-1)^{(k-1)p}}{|H|} \sum_{\si\in\si_0 H} \chi^{\la}(\si),
$$
\end{lem}
\begin{proof}
Write 
$$
\{1,\ldots,kp\} = C_1 \cup \ldots \cup C_k,
$$
where
$$
C_i := \{(i-1)p+1, (i-1)p+2, \ldots, (i-1)p+p\} \qquad(1\le i\le k).
$$
This decomposition has the property that $j\in C_i$ if and only if the $j$th
factor in the tensor product (\ref{e.v0}) is $e_i$.

Now clearly, for any $\si\in S_{kp}$
$$
\langle \htau\rho(\si)(v_0), v_0 \rangle \in \{0,(-1)^{(k-1)p}\},
$$
where the sign comes from the action of $\htau$ on $p$ copies of $e_k$.
Indeed, by our choice of $\htau$, 
$$
\htau(e_i) = \pm e_{i+1} \qquad(1\le i\le k)
$$
where indices (here and in the sequel) are computed modulo $k$. Thus
$\langle \htau\rho(\si)(v_0), v_0 \rangle \ne 0$ if and only if 
$\htau\rho(\si)(v_0) = \pm v_0$, that is 
$$
\si(C_i) = C_{i+1} \qquad(1\le i\le k).
$$
This happens if and only if there are $k$ permutations 
$\si_1,\ldots,\si_k\in S_p$ such that
$$
\si((i-1)p+j) = ip+\si_{i}(j) \qquad(1\le i\le k,\, 1\le j\le p).
$$
Let $\si_0 \in S_m$ correspond to $\si_1 = \ldots = \si_k = id \in S_p$, 
so that
$$
\si_0((i-1)p+j) = ip+j \qquad(1\le i\le k,\, 1\le j\le p).
$$
Then
$$
\langle \htau\rho(\si)(v_0),v_0 \rangle \ne 0 \iff
\si_0^{-1} \si(C_i) = C_i\quad(\forall i) \iff 
\si_0^{-1} \si \in H,
$$
and we conclude from Lemma~\ref{t.eps.sm} that
$$
\trt = \frac{(-1)^{(k-1)p}}{|H|} \sum_{\si_0^{-1}\si\in H} \chi^{\la}(\si),
$$
as claimed.
\end{proof}

\subsection{From $S_p^{\times k}$ to $S_p$}

Denote now
$$
\epsilon_1(\la) := \frac{1}{|H|} \sum_{\si\in\si_0 H} \chi^\la(\si),
$$  
so that
$$
\trt = (-1)^{(k-1)p} \epsilon_1(\la).
$$
Recall that each dominant weight $\bla$ may be obtained from many partitions 
$\la\in Par_k(m)$.  There is a minimal such $\la = (\la_1,\ldots,\la_k)$, 
with $\la_k = 0$.  All the other partitions are 
obtained from it by adding the same integer $c$ to each of 
$\la_1,\ldots,\la_k$, thus adding $c$ to $p = (\la_1 + \ldots + \la_k)/k$.

\begin{lem}\label{t.trt.chi}
If $\la\vdash kp$ then
$$
\epsilon_1(\la) = \sum_{\mu\vdash p} \frac{1}{z_{\mu}} \chi_{k\mu}^{\la}.
$$
In particular, if $\core_k(\la) \ne \emptyset$ then $\epsilon_1(\la) = 0$.
\end{lem}
\begin{proof}
Recall the definition of $\si_0 \in S_{kp}$ from the proof of Lemma~\ref{t.sm.h}.
For each $\si\in\si_0 H$, $\si(C_i) = C_{i+1}$ ($1\le i\le k$) so that
all the cycle lengths of $\si$ are divisible by $k$.
As in the proof of Lemma~\ref{t.sm.h}, there exist $k$ permutations
$\si_1,\ldots,\si_k\in S_p$ such that
$$
\si((i-1)p+j) = ip+\si_{i}(j) \qquad(1\le i\le k,\, 1\le j\le p).
$$
The character value $\chi^{\la}(\si)$ depends only on the cycle lengths
of $\si$, which in turn depend only on the product 
$\bsi:= \si_k \cdots \si_1 \in S_p$.  To each cycle of length $t$ in $\bsi$ 
there corresponds a cycle of length $kt$ in $\si$.

Define a map $\phi: S_p^{\times k} \too S_p$ by 
$$
\phi(\si_1,\ldots,\si_k) := \si_k \cdots \si_1\qquad(\forall \si_1,\ldots,\si_k\in S_p).
$$
Clearly, the size of the inverse image
$$
|\phi^{-1}(\bsi)| = (p!)^{k-1} = \frac{|H|}{|S_p|}\qquad(\forall \bsi\in S_p). 
$$
If $\bsi\in S_p$ has cycle type $\mu$, then each $\si\in\si_0 H$ corresponding
to a member of $\phi^{-1}(\bsi)$ has cycle type $k\mu$, so that
$\chi^\la(\si) = \chi_{k\mu}^\la$. 
By Proposition~\ref{t.cycletype}, there are $|S_p|/z_{\mu}$ permutations in $S_p$
of any given cycle type $\mu$, so that
$$
\epsilon_1(\la) := \frac{1}{|H|} \cdot  \frac{|H|}{|S_p|} 
\sum_{\mu\vdash p} \frac{|S_p|}{z_{\mu}} \chi_{k\mu}^\la
= \sum_{\mu\vdash p} \frac{1}{z_{\mu}} \chi_{k\mu}^\la,
$$
as claimed.
\end{proof}

\begin{thm}\label{t.eps.sign}\
\begin{description}
\item{(a)}
If $\la\vdash kp$ has $\core_k(\la) = \emptyset$ and $\quot_k(\la) = \lala$, 
where each partition $\la_i$ is either empty or has one part, then
$$
\epsilon_1(\la) = \sign(\pi_{\la})
$$
and
$$
\trt = (-1)^{(k-1)p} \sign(\pi_\la).
$$
Here $\pi_\la$ is the zero permutation of $\la$, as in Lemma~\ref{t.zero_perm}.
\item{(b)}
For any other partition $\la$,
$$
\trt = \epsilon_1(\la) = 0.
$$  
\end{description}
\end{thm}
\begin{proof}
Assume that $\la\vdash kp$ and $\core_k(\la) = \emptyset$; otherwise we are in 
case (b) of our theorem, and $\epsilon_1(\la) = 0$ 
by Lemma~\ref{t.basis}(a), equation~(\ref{e.trt}), and Lemma~\ref{t.trt.chi}.
Now, by Lemma~\ref{t.trt.chi}, Proposition~\ref{t.mnf} (the Murnaghan-Nakayama 
formula), and Lemma~\ref{t.sign.hgt}:
\begin{eqnarray*}
\epsilon_1(\la) 
     &=& \sum_{\mu\vdash p} \frac{1}{z_{\mu}} \chi_{k\mu}^{\la} =\\
     &=& \sum_{\mu\vdash p} \frac{1}{z_{\mu}} 
         \sum_{T\in\rht_{k\mu}^{\la}} (-1)^{\hgt(T)} =\\
     &=& \sum_{\mu\vdash p} \frac{1}{z_{\mu}} 
         \sum_{T\in\rht_{k\mu}^{\la}} \sign(\pi_T).
\end{eqnarray*}
By Lemma~\ref{t.size}, 
$\quot_k(\la) = \lala$ satisfies
$$
|\la_1| + \ldots + |\la_k| = p \quad(= |\la|/k).
$$ 

Denote now $p_i := |\la_i|$ $(1\le i\le k)$.
Recalling the definition of $\omega_{k\mu}^{\la}$ from Subsection~\ref{s.quot.tab},
let $\omega_{k\mu}^{\la}(T) = (\TT)$.
By Lemma~\ref{t.rht.dec} and Theorem~\ref{t.sign}:
$$
\epsilon_1(\la)
       = \sum_{\mu\vdash p} \frac{1}{z_{\mu}} \sum_{\mumu} 
         \frac{z_{\mu}}{z_{\mu_1}\cdots z_{\mu_k}} 
         \sum_{\TT} \sign(\pi_{\la}) \prod_{i=1}^{k} \sign(\pi_{T_i}),
$$
where the second summation is over all $\mu_i\vdash p_i$ such that 
$\mu_1 \oplus \ldots \oplus \mu_k = \mu$ and the third summation is
over all $T_i \in \rht_{\mu_i}^{\la_i}$.  Equivalently,
$$
\epsilon_1(\la)
     = \sign(\pi_{\la}) \cdot \sum_{\mumu} \prod_{i=1}^{k} \frac{1}{z_{\mu_i}} 
       \sum_{T_i \in \rht_{\mu_i}^{\la_i}} \sign(\pi_{T_i}),
$$
where the first summation is over all $\mu_i \vdash p_i$.
Thus, by Lemma~\ref{t.sign.hgt} and 
Propositions~\ref{t.mnf} and~\ref{t.charsum}:
\begin{eqnarray*}
\epsilon_1(\la)
     &=& \sign(\pi_{\la}) \cdot \prod_{i=1}^{k} 
         \sum_{\mu_i\vdash p_i} \frac{1}{z_{\mu_i}} \chi_{\mu_i}^{\la_i} =\\
     &=& \sign(\pi_{\la}) \cdot \prod_{i=1}^{k} \delta_{\la_i,(p_i)}.
\end{eqnarray*}
This completes the proof.
\end{proof}

This result can be given a very explicit reformulation.

\begin{lem}
For a partition $\la \vdash kp$, $\epsilon_1(\la)\ne 0$ if and only if $\la$
has at most $k$ nonzero parts: $\la_1 \ge \ldots \ge \la_k \ge 0$, and the 
numbers $(\la_i + k - i)_{i=1}^{k}$ have $k$ distinct residues $(\mmod k)$.
In this case, $\epsilon_1(\la)$ is the sign of the permutation needed to 
transform the sequence of residues $(\mmod k)$ of $(\la_i + k - i)_{i=1}^{k}$ 
into the sequence $(k - i)_{i=1}^{k}$.  
\end{lem}
\begin{proof}
By Subsection~\ref{s.01.sub}, a partition $\la$ has at most one part 
if and only if, in its $0/1$ sequence $s(\la)$, all the digits before the median,
except possibly the last one, are zeros. Thus, by Lemma~\ref{t.empty.core},
$\la$ satisfies the assumptions of Theorem~\ref{t.eps.sign}(a) if and only if:
$s(\la)$ has only zeros until the $k$th ``space'' before the median; 
has exactly $k$ zeros after this ``space $-k$''; 
and the indices in $s(\la)$ of these $k$ zeros have distinct residues 
$(\mmod k)$. 
In particular, $\la$ has at most $k$ parts.
The index in $s(\la)$ of the zero corresponding to the $i$th part of $\la$
$(1\le i\le k)$ is $\la_i+k-i$, where index $0$ corresponds to the digit 
immediately following ``space $-k$''; in $s(\emptyset)$, the index is $k-i$.
Thus the zero permutation $\pi_{\la}$ is exactly the permutation described in
the statement of the lemma.
\end{proof}

Finally, recall that $\trt = (-1)^{(k-1)p}\epsilon_1(\la)$. Adding $1(\mmod k)$ to
each of $k$ distinct residues is equivalent to permuting them by a cycle of 
length $k$, whose sign is $(-1)^{k-1}$. Thus we get the following result.

\begin{cor}\label{t.bla.final}
Let $\bla = (\bla_1,\ldots,\bla_{k-1})$ be a dominant $SU(k)$-weight, and define
$\bla_k := 0$. If the numbers $(\bla_i + k - i)_{i=1}^{k}$ 
have $k$ distinct residues $(\mmod k)$
then $\trt$ is the sign of the permutation needed to transform the sequence 
of residues $(\mmod k)$ of $(\bla_i -\bar{p} + k - i)_{i=1}^{k}$ into 
the sequence $(k-i)_{i=1}^{k}$, where $\bar{p} := (\bla_1 + \ldots + \bla_k)/k$.
In all other cases, $\trt = 0$.
\end{cor}

{\bf Acknowledgments.} The authors thank Bertram Kostant and Amitai Regev 
for useful comments.

\end{document}